\documentclass{article}%
\usepackage{amssymb}
\usepackage{latexsym}
\usepackage{amsmath}
\usepackage{graphicx}
\usepackage{amsfonts}%
\newtheorem{theorem}{Theorem}
\newtheorem{acknowledgement}[theorem]{Acknowledgement}

\newtheorem{conjecture}[theorem]{Conjecture}
\newtheorem{corollary}[theorem]{Corollary}
\newtheorem{criterion}[theorem]{Criterion}
\newtheorem{definition}[theorem]{Definition}

\newtheorem{lemma}[theorem]{Lemma}
\newtheorem{notation}[theorem]{Notation}

\newtheorem{proposition}[theorem]{Proposition}
\newtheorem{remark}[theorem]{Remark}

\begin{document}

\author{Andrey N. Todorov\\UC, Department of Mathematics\\Santa Cruz, CA 95064,USA\\Bulgarian Academy of Sciences\\Institute of Mathematics\\Sofia, Bulgaria}
\title{Large Radius Limit and SYZ Fibrations of Hyper-K\"{a}hler Manifold}
\date{}
\maketitle

\begin{abstract}
In this paper the relations between the existence of Lagrangian fibration of
Hyper-K\"{a}hler manifolds and the existence of the Large Radius Limit is
established. It is proved that if the the rank of the second homology group of
a Hyper-K\"{a}hler manifold N of complex dimension $2n\geq4$ is at least 5,
then there exists an unipotent element T in the mapping class group $\Gamma
$(N) such that its action on the second cohomology group satisfies
$(T-id)^{2}\neq0$ and $(T-id)^{3}=0.$ A Theorem of Verbitsky implies that the
symmetric power $S^{n}(T)$ acts on $H^{2n}$ and it satisfies $(S^{n}%
(T)-id)^{2n}\neq0$ and $(S^{n}(T)-id)^{2n+1}=0.$ This fact established the
existence of Large Radius Limit for Hyper-K\"{a}hler manifolds for polarized
algebraic Hyper-K\"{a}hler manifolds. Using the theory of vanishing cycles it
is proved that if a Hyper-K\"{a}hler manifold admits a Lagrangian fibration
then the rank of the second homology group is greater than or equal to five.
It is also proved that the fibre of any Lagrangian fibration of a
Hyper-K\"{a}hler manifold is homological to a vanishing invariant $2n$ cycle
of a maximal unipotent element acting on the middle homology. According to
Clemens this vanishing invariant cycle can be realized as a torus. I
conjecture that the SYZ conjecture implies finiteness of the topological types
of Hyper-K\"{a}hler manifolds of fix dimension.

\end{abstract}
\tableofcontents

\section{Introduction}

The mirror symmetry conjecture attracted the attention of mathematicians after
the appearance of the paper by Ph. Candelas and his coauthors. See \cite{Can}.
The mirror symmetry conjecture is based on the suggestion that two
sigma-models in superstring theory are equivalent. The targets of the models
are Calabi-Yau 3-folds, so mirror symmetry predicts that these should come in
pairs, M and \v{M}, satisfying $h^{p,q}($M)=$h^{p,3-q}$(\v{M}). Recently
Strominger-Yau-Zaslow noticed in \cite{SYZ} that String Theory suggests that
if M and \v{M} are mirror pairs of n-dimensional CY manifolds, then on M
should exist a special Lagrangian n-tori fibration f:M$\rightarrow$B, (with
some singular fibres) such that \v{M} is obtained by finding some suitable
compactification of the dual fibration. That fibration, which I will call the
SYZ fibration, should exist near the large radius limit points in the moduli
space of Calabi-Yau manifolds. Recently Mark Gross and Ruan have made some
progress with respect to the topological properties of the SYZ fibration of
the quintic CY in $\mathbb{CP}^{4}.$ See \cite{Gr} and \cite{R}. For K3
surfaces a proof of the existence of SYZ fibrations was given by M. Gross.

Some aspects of mirror symmetry for Hyper-K\"{a}hler manifolds were developed
by M. Verbitsky in \cite{Ver}. Verbitsky showed that a compact
Hyper-K\"{a}hler manifold is mirror to itself in a strong sense. He showed
that the Yukawa coupling is just the ordinary cup product.

I am not planning to discuss the Fourier-Mukai transforms of supersymmetric
cycles, which is part of SYZ approach to mirror symmetry. Some very
interesting results about these aspects of mirror symmetry were obtained in
\cite{LYZ} and \cite{ACRY}. For detailed account of the recent developments in
Mirror Symmetry see \cite{L1}, \cite{L2} and \cite{Man}. Some basic facts
about Lagrangian fibrations are discussed in the papers \cite{H1}, \cite{H2}
and \cite{H3}. A very interesting discussion about infinite dimensional point
of view of SYZ fibrations can be found ins the paper of S. Donaldson
\cite{Don}.

Not too much is known about the existence of such fibrations. In the important
paper \cite{HT}, B. Hassett and Yu. Tschinkel discuss the existence of abelian
fibrations on holomorphic symplectic fourthfolds with Picard group of rank 2
and the authors gave examples of such fibrations on fourthfolds.

In this paper some aspects of SYZ conjecture about compact Hyper-K\"{a}hler
manifold will be discussed. Hyper-K\"{a}hler manifold is defined as follows:

\begin{definition}
Let N be a compact K\"{a}hler manifold of complex dimension $2n\geq4$ such
that on N there exists a non zero holomorphic two form $\Omega_{\text{N}},$
\ unique up to a constant such that $\det(\Omega_{\text{N}})\neq0$ at each
point $m\in$N and
\[
H^{1}(\text{N},\mathcal{O}_{\text{N}})=0.
\]
Then N will be called a Hyper-K\"{a}hler manifold.
\end{definition}

The conditions on the holomorphic two form $\Omega_{\text{N}}$ imply that
\[
\dim_{\mathbb{C}}H^{2}(\text{N,}\mathcal{O}_{\text{N}})=1.
\]

In the case of Hyper-K\"{a}hler manifold N of complex dimension $2n,$ part of
the conjecture of Strominger-Yau-Zaslow states that if one deforms the complex
structure isometrically with respect to some Calabi Yau metric then a new
complex Hyper-K\"{a}hler structure N$_{0}$ on N is obtained. Moreover N$_{0}$
admits a fibration over $\mathbb{CP}^{n}$ by complex Lagrangian n-dimensional
tori $\mathbb{T}^{n}.$ Such a fibration will be called SYZ fibration. It is
suggested in the Physics literature that the problem of the existence SYZ
fibrations is related to the existence of large radius limit points in the
moduli space of CY manifolds. The main result of this article is to establish
the relation between the existence of the large radius limit and the existence
of the Lagrangian fibration on Hyper-K\"{a}hlerian manifold N.

We will use the following definition for the large limit radius:

\begin{definition}
\label{LRL}We will say the moduli space of polarized algebraic
Hyper-K\"{a}hlerian manifold N admits a large radius limit if there exists an
element $T$ of the mapping class group $\Gamma_{\text{N}}$ such that the
action of $T$ on the middle cohomology $H^{2n}($N,$\mathbb{Z})$ satisfies:
\begin{equation}
(T-id)^{2n+1}=0\text{ and }(T-id)^{2n}\neq0. \label{u2}%
\end{equation}

\end{definition}

We will prove that Definition \ref{LRL} is equivalent to the following definition:

\begin{definition}
\label{LRLI}There exists a family of Hyper-K\"{a}hler polarized algebraic
manifolds%
\[
\mathcal{X\rightarrow D}%
\]
over the unit disc such that the action of the monodromy operator of this
family $T$ satisfies $\left(  \ref{u2}\right)  $.
\end{definition}

In this paper it is proved the existence of large radius limit if
\begin{equation}
rkH^{2}(\text{N,}\mathbb{Z})\geq5. \label{u1}%
\end{equation}
We also established that if the Hyper-K\"{a}hlerian manifold N admits a
Lagrangian fibration than $\left(  \ref{u1}\right)  $ holds.

One of the observation made in this paper is that the existence of large
radius limit of a Hyper-K\"{a}hlerian manifold N is equivalent to the
existence of an isotropic non zero vector $\delta\in H^{2}($N,$\mathbb{Z})$
with respect to Beauville-Bogomolov form defined in Definition \ref{Bog-B}.
Then it is not difficult to see that the subgroup of the automorphisms of the
Euclidean lattice defined by $H^{2}($N,$\mathbb{Z})/Tor$ and
Bogomolov-Beauville quadratic form that stabilizes the isotropic vector
$\delta$ will contain an automorphism $T_{1}$ such that it satisfies%
\[
(T_{1}^{m}-id)^{2}\neq0\text{ and }(T_{1}^{m}-id)^{3}=0.
\]
See \cite{Sc}. Then a Theorem of Verbitski combined with some elementary
linear algebra imply that the $n$ symmetric power $T:=S^{n}T_{1}$ acts on
$H^{2n}($N,$\mathbb{Z})$ and satisfies $\left(  \ref{u2}\right)  $. A relation
is established between the monodromy operator $T$ of maximal index of
unipotency and the Lagrangian fibration by using Clemens theory of vanishing
cycles. See \cite{LTY}. It is proved that the homology class $[\mathbb{T}%
^{n}]$ of the Lagrangian fibration is the same as the homology class of the
vanishing invariant cycle of the monodromy operator $T^{m}.$ It seems
reasonable to expect that the existence of the large radius limit and Clemens
theory will imply the existence of a topological fibration of Hyper-K\"{a}hler
manifolds by real tori of half of the dimension of N.

It is proved that if $\mathcal{L}_{l}$ is the line bundle with Chern class
$c_{1}(\mathcal{L}_{l})=l\in H^{2}($N,$\mathbb{Z})$ and it is isotropic with
respect to the Bogomolov-Beauville quadratic form then the Euler
characteristics $\chi(\mathcal{L}_{l})$ is equal to $n+1,$ i.e.
\[
\chi(\mathcal{L}_{l})=%
{\displaystyle\sum\limits_{i=0}^{2n}}
\dim H^{i}(\text{N,}\mathcal{L}_{l})=n+1.
\]
This result suggests that if N is a hyper K\"{a}hler manifold such that its
Picard group $H^{2}($N,$\mathbb{Z})\cap H^{1,1}($N,$\mathbb{R})$ has rank one,
then probably the linear system \ defined by $H^{0}($N,$\mathcal{L}_{l})$ will
define the SYZ fibration. It is also can be established that if $\delta\in
H^{2}($N,$\mathbb{Z})$ is a non zero isotropic vector of type $(1,1)$ and the
linear system \ defined by $H^{0}($N,$\mathcal{L}_{l})$ defines the SYZ
fibration then the Poincare dual class of homology $\mathcal{P}(\wedge
^{n}\delta$ ) of $\wedge^{n}\delta$ will coincide with the class of homology
of the fibre of the SYZ fibration.

I believe that the construction of the SYZ fibration can be used to classify
the compact Hyper-K\"{a}hler manifolds or to try to prove that there are only
a finite number of topological types of Hyper-K\"{a}hler manifolds in
dimension 4. It seems that the existence of SYZ fibrations can lead to the
topological classification of Hyper-K\"{a}hler manifold in particular in
dimension 4, since it implies that any Hyper-K\"{a}hler manifold in dimension
4 can be isometrically deformed to a Hyper-K\"{a}hler manifold which is
fibered by tori over $\mathbb{CP}^{2}$. It is reasonable to expect that one
can find a bound for the degree of the divisor of singular fibres. Once such
bound is established then I expect that the classification all such fibrations
will not be so difficult.

\begin{acknowledgement}
This paper was written during my stay at the Mathematical Institute of the
Chinese Academy of Sciences in Beijing and my visit to the Institute of
Mathematical Sciences at CUHK. I want to thank my Chinese colleagues for their
hospitality and excellent working conditions. I want to thank Professor Yau
for his generous advice, encouragements and help. Special thanks to F.
Bogomolov, Yu. Tchinkel, B. Hassett and C. Leung for their help. The comments
of the referee were also very useful and helpful. I would like to thank Jamey
Bass for his help in the preparation of the article.
\end{acknowledgement}

\section{Basic Definitions}

In this paper I will assume that the Hyper-K\"{a}hler manifolds have complex
dimension $\geq4.$

\begin{definition}
\label{d2}Let $\gamma_{1},...,\gamma_{b_{2}}$ be a basis in $H_{2}%
($N,$\mathbb{Z})/Tor,$ then the pair%

\[
(\text{N};\gamma_{1},...,\gamma_{b_{2}})
\]
will be called a marked Hyper-K\"{a}hler manifold. Suppose that $L\in H^{2}%
($N$,\mathbb{Z})$ can be realized as an imaginary part of a K\"{a}hler metric
on N. Then the triple (N;$L$;$\gamma_{1},...,\gamma_{b_{2}}$) will be called a
marked polarized Hyper-K\"{a}hler manifold.
\end{definition}

\begin{definition}
\label{q}For every pair of cocycles $\alpha$ and $\beta\in H^{2}%
($N$,\mathbb{R)}$ I define
\begin{equation}
q(\alpha,\beta):=\int_{\text{N}}\alpha\wedge\beta\wedge L^{2n-2}. \label{-q1}%
\end{equation}

\end{definition}

\begin{proposition}
\label{q1}The quadratic form $q$ is a non-degenerate one. It has a signature
$(2,b_{2}-3)$ on the primitive cohomology
\begin{equation}
H_{0}^{2}(\text{N,}\mathbb{Z}):=\left\{  \alpha\in H^{2}(\text{N,}%
\mathbb{Z})|\text{ }q(\alpha,L)=0\right\}  \label{prim}%
\end{equation}
and it is defined over $\mathbb{Z}.$
\end{proposition}

\textbf{Proof: }The Proposition \ref{q1} is an easy exercise. $\blacksquare$

\begin{notation}
If N is a Hyper-K\"{a}hler manifold, then I will denote by $\Omega_{\text{N}}$
the non-zero holomorphic two form. If $\pi:\mathcal{N\rightarrow C}$ is a
family of Hyper-K\"{a}hler manifolds and if $\tau\in\mathcal{C}$ is a point on
the base, I will denote by $\Omega_{\tau}$ the non zero holomorphic two form
on the fibre N$_{\tau}=\pi^{-1}(\tau).$
\end{notation}

\begin{definition}
\label{Bog-B}For every $\alpha\in H^{2}($N$,\mathbb{R)}$ let us define
\[
\mathcal{B}(\alpha):=\frac{n}{2}\int_{\text{N}}\left(  \wedge^{n-1}%
\Omega_{\text{N}}\right)  \wedge\left(  \wedge^{n-1}\overline{\Omega
_{\text{N}}}\right)  \wedge\alpha^{2}+
\]%
\begin{equation}
(1-n)\left(  \int_{\text{N}}\left(  \wedge^{n-1}\Omega_{\text{N}}\right)
\wedge\left(  \wedge^{n}\overline{\Omega_{\text{N}}}\right)  \wedge
\alpha\right)  \left(  \int_{\text{N}}\left(  \wedge^{n}\Omega_{\text{N}%
}\right)  \wedge\left(  \wedge^{n-1}\overline{\Omega_{\text{N}}}\right)
\wedge\alpha\right)  . \label{q2}%
\end{equation}
I will call $\mathcal{B}(\alpha)$ Beauville-Bogomolov form.
\end{definition}

\begin{proposition}
\label{q3}The Beauville-Bogomolov form $\mathcal{B}$ is a non-degenerate one.
It has a signature $(3,b_{2}-3)$ and it is defined over $\mathbb{Z}$ after
multiplying it with a suitable constant.
\end{proposition}

\textbf{Proof: }The Proposition \ref{q3} was proved by Beauville. See
\cite{Bau}. $\blacksquare$

\begin{remark}
\label{BB0}The forms defined in $\left(  \ref{-q1}\right)  $ and in $\left(
\ref{q2}\right)  $ coincide on the primitive cohomology classes of type
$(1,1)$ on N with respect to the polarization class $L$ up to a constant$.$
\end{remark}

\begin{notation}
I will denote the intersection of any two cycles $\alpha$ and $\beta\in
H_{\ast}($N,$\mathbb{Z})$ by $\left\langle \alpha,\beta\right\rangle .$ If
$\alpha$ has dimension $p$ and $\beta$ has a dimension $q$ and $p+q\geq4n$,
then $\left\langle \alpha,\beta\right\rangle $ will have dimension $p+q-4n.$
\end{notation}

\section{The Existence of Large Radius Limit for Hyper-K\"{a}hler Manifolds}

\subsection{Review of Some Elementary Facts about Unipotent Elements in the
Mapping Class Group}

\begin{definition}
\label{MSG}We will define the mapping class group $\Gamma$(N) of N as follows%
\[
Diff^{+}\text{(N)/}Diff_{0}(\text{N)}%
\]
where $Diff^{+}$(N) is the group of diffeomorphisms of N that preserve the
orientation and $Diff_{0}$(N) is the group of diffeomorphisms of N isotopic to
identity. The subgroup $\Gamma_{L}$ of $\Gamma($N) will be defined as
\[
\Gamma_{L}:=\left\{  \phi\in\Gamma\text{(N)%
$\vert$%
}\phi(L)=L\right\}
\]
where $L$ is a fixed class in $H^{2}($N,$\mathbb{Z}).$
\end{definition}

According to the results of Sullivan both groups $\Gamma$(N) and $\Gamma_{L}$
are arithmetic subgroups mapped to $SO(3,n)$ and $SO(2,n)$ with finite kernels
respectively$,$ where
\[
n+3=rkH^{2}(\text{N},\mathbb{Z}).
\]
Next we will prove some elementary facts about the relations between the
Jordan blocks of a linear operator acting on a vector space with the Jordan
blocks of its symmetric power.

\begin{theorem}
\label{Mon1}Suppose that $\phi\in\Gamma$(N) and $\phi$ has a maximal index of
unipotency, then the induced action $T_{\phi}$ of $\phi$ on $H^{2}%
($N,$\mathbb{Z})$ has only one Jordan block of dimension three and so it
satisfies:
\begin{equation}
(T_{1}-id)^{2}\neq0\text{ and }(T_{1}-id)^{3}=0. \label{Uni}%
\end{equation}
The $n$ symmetric power $T=S^{n}(T_{\phi})$ of $T_{\phi}$ has only one Jordan
block of dimension $2n+1$ on $H^{2n}($M,$\mathbb{Z})$ and so $T$ has the
following properties$:$
\[
(T-id)^{2n}\neq0\text{ and }(T-id)^{2n+1}=0.
\]

\end{theorem}

\textbf{Proof: }The proof of Theorem \ref{Mon1} is based on the following
important Theorem due to Verbitsky:

\begin{theorem}
\label{Ver}Let N be Hyper-K\"{a}hler manifold of dimension $2n$. Then
\begin{equation}
S^{k}H^{2}(\text{N})\subseteq H^{2k}(\text{N}) \label{V1}%
\end{equation}
for $k\leq n.$ Bogomolov found a very elegant proof of Theorem \ref{Ver} in
\cite{Bo1}.
\end{theorem}

and the following remark:

\begin{remark}
\label{lty}In \cite{LTY} \ it was proved that if $\phi\in\Gamma$(N) and the
induced action of $T_{\phi}$ of $\phi$ on $H^{2k}($N,$\mathbb{Z})$ for $k\leq
n$ has a maximal index of unipotency then $T_{\phi}$ has only one Jordan block
of dimension $2k+1.$
\end{remark}

\begin{lemma}
\label{031}Suppose that $V$ is a vector space over $\mathbb{C}$ and $T_{\phi}$
is a linear operator acting on $V$ such that it has only one Jordan block of
dimension $3.$ Then the n$^{th}$ symmetric power $T=S^{n}(T_{\phi})$ of the
operator $T_{\phi}$ will have only one Jordan block of dimension 2n+1.
\end{lemma}

\textbf{Proof: }The assumption that the linear operator $T_{\phi}$ contains
only one Jordan block of dimension 3 implies that there exist non zero vectors
$v_{0},v_{1}$ and $v_{2}$ such that
\begin{equation}
T_{\phi}(v_{0})=v_{0},T_{\phi}(v_{1})=v_{0}+v_{1}\text{ and }T_{\phi}%
(v_{2})=v_{2}+v_{1}. \label{Jor3}%
\end{equation}
Then $\left(  \ref{Jor3}\right)  $ implies that the following $2n+1$ vectors
in $V:$%
\[
\gamma_{0}=\underset{n}{\underbrace{\left\langle v_{0},...,v_{0}\right\rangle
}},\text{ }\gamma_{1}=\underset{n}{\underbrace{\left\langle v_{0}%
,...,v_{0},v_{1}\right\rangle }},...,\text{ }\gamma_{n}=\underset
{n}{\underbrace{\left\langle v_{0},v_{1},...,v_{1}\right\rangle }},
\]%
\begin{equation}
\gamma_{n+1}=\underset{n}{\underbrace{\left\langle v_{1},v_{1},...,v_{1}%
\right\rangle }},\text{ }\gamma_{n+2}=\underset{n}{\underbrace{\left\langle
v_{1},...,v_{1},v_{2}\right\rangle }},...,\text{ }\gamma_{2n}=\underset
{n}{\underbrace{\left\langle v_{2},...,v_{2}\right\rangle }} \label{M0}%
\end{equation}
are linearly independent. Direct computation shows that the symmetric power
$S^{n}(T_{\phi})=T$ acts on $\gamma_{i}$ for $i=0,...,2n$ as follows:
\begin{equation}
T(\gamma_{0})=\gamma_{0},...,\text{ }T(\gamma_{i})=\gamma_{i}+\gamma
_{i-1},...,\text{ }T(\gamma_{2n})=\gamma_{2n}+\gamma_{2n-1}. \label{Jor1}%
\end{equation}
Lemma \ref{031} follows directly from $\left(  \ref{Jor1}\right)  $.
$\blacksquare$

Theorem \ref{Mon1} follows directly from Lemma \ref{031}, Remark \ref{lty} and
Theorem \ref{Ver}. $\blacksquare$

\subsection{Existence of Large Radius Limit}

\begin{theorem}
\label{Cond}Suppose that N is a Hyper-K\"{a}hler manifold such that%
\begin{equation}
\dim_{\mathbb{Q}}H^{2}(N,\mathbb{Q})\geq5. \label{Cond1}%
\end{equation}
Then the large limit radius exists for N in the since of Definition \ref{LRL}.
\end{theorem}

\textbf{Proof:} The proof of Theorem \ref{Cond} is based on the following
Theorem of Meyer:

\begin{theorem}
\label{Me} Let $\Lambda$ be an Euclidean lattice whose scalar product is
defined by an indefinite quadratic form. Suppose that
\[
rk_{\mathbb{Z}}\Lambda\geq5.
\]
Then there exists a non zero vector $\delta\in\Lambda$ such that
\[
\left\langle \delta,\delta\right\rangle =0.
\]
(For the proof of this Theorem see \cite{Me} or \cite{MH}.)
\end{theorem}

The condition $\left(  \ref{Cond1}\right)  $ implies that we can apply Theorem
\ref{Me} to deduce that there exists a non zero vector%
\[
\delta\in H^{2}\text{(N,}\mathbb{Z})
\]
which is isotropic with respect to the Beauville-Bogomolov form, i.e.
\begin{equation}
\mathcal{B}(\delta,\delta)=0. \label{Cond2}%
\end{equation}
From the results proved in \cite{Sc} and the fact that the mapping class group
$\Gamma$(N) is an arithmetic group we know that the subgroup of the mapping
class group $\Gamma$(N) whose intersection with the automorphism group of the
Euclidean lattice $H^{2}$(N,$\mathbb{Z})$ with respect to the
Beauville-Bogomolov form has finite indices in both groups contains an
unipotent element $T_{\phi}$ which has a unique Jordan block of dimension
three. The uniqueness of only one the Jordan block of $T_{\phi}$ follows from
Remark \ref{lty}. Thus the action of $T_{\phi}$ on the lattice $H^{2}%
$(N,$\mathbb{Z})$ satisfies the conditions of Theorem \ref{Mon1}. Theorem
\ref{Cond} is proved. $\blacksquare$

\subsection{Relations between the Lagrangian Fibered Hyper-K\"{a}hler
Manifolds and Large Radius Limit}

\begin{theorem}
\label{Top}Suppose that N is a Hyper-K\"{a}hler Manifold that admits a
Lagrangian fibration
\begin{equation}
\pi:\text{N}\rightarrow\text{B,} \label{Lagfibr}%
\end{equation}
where B is a complex manifold. Then the base\ B of the fibration is a
projective variety and
\[
\dim_{\mathbb{Q}}H^{2}(\text{N},\mathbb{Q})\geq5.
\]

\end{theorem}

\textbf{Proof:} The proof of Theorem \ref{Top} will be based on the following Lemma:

\begin{lemma}
\label{ac1}The base B of the fibration $\left(  \ref{AC0}\right)  $ is an
algebraic manifold.
\end{lemma}

\textbf{Proof: }Since N is a K\"{a}hler manifold, I can conclude that the
image $\pi($N)$=$B is also K\"{a}hler by a well known result due to Varochas.
On the other hand it is easy to see that
\begin{equation}
H^{0}(\text{B},\Omega_{\text{B}}^{2})=0. \label{ac10}%
\end{equation}
Indeed since B is a complex manifold and if on B there exists a holomorphic
two form $\omega$, then $\pi^{\ast}(\omega)$ will be a non zero holomorphic
form on N. From the fact that the holomorphic two form is lifted from the base
I can conclude that $\pi^{\ast}(\omega)$ is a degenerate two form on N. Indeed
let $x\in$N be a point in some non-singular fibre $\pi^{-1}(t)=$
$\mathbb{T}_{t}.$ Let $u,v\in T_{x,\text{N}}$ , where $u$ is tangent to the
fibre $\pi^{-1}(t)$ and $v$ is such that
\begin{equation}
\Omega_{\text{N}}(u,v)\neq0. \label{def0}%
\end{equation}
The existence of $u$ and $v$ follows from the fact that we assumed that the
fibration is Lagrangian and $\Omega_{\text{N}}$ is the non degenerate
holomorphic two form on N. From the way the holomorphic two form $\pi^{\ast
}(\omega)$ was defined it follows that
\begin{equation}
\pi^{\ast}(\omega)\left\langle u,v\right\rangle =0. \label{deg}%
\end{equation}
The condition $\left(  \ref{deg}\right)  $ implies that the holomorphic two
form $\pi^{\ast}(\omega)$ is degenerate on N, i.e.
\begin{equation}
\det(\pi^{\ast}(\omega))=0. \label{deg-0}%
\end{equation}
Since on N there exists a holomorphic two form $\Omega_{\text{N}}$ unique up
to a constant which is a non degenerate one I can conclude from $\left(
\ref{deg-0}\right)  $ and the fact that N is a K\"{a}hler manifold that
\begin{equation}
\dim_{\mathbb{C}}H^{0}(\text{N},\Omega_{\text{N}}^{2})\geq2. \label{ac11}%
\end{equation}
But $\left(  \ref{ac11}\right)  $ contradicts the definition of a
Hyper-K\"{a}hler manifold, since it states that
\begin{equation}
\dim_{\mathbb{C}}H^{0}(\text{N},\Omega_{\text{N}}^{2})=1. \label{ac12}%
\end{equation}
The condition $\left(  \ref{ac12}\right)  $ implies that B is a K\"{a}hler
manifold with $H^{2,0}($B$)=0.$ This fact and the fact that the K\"{a}hler
condition is an open one imply that I can find a class of cohomology of type
$(1,1)$ with integral periods which can be realized as the imaginary class of
a K\"{a}hler metric. Kodaira proved that such manifolds are projective
algebraic. See \cite{KM}. Lemma \ref{ac1} is proved. $\blacksquare$

\begin{lemma}
\label{ac-1}Let $D_{\text{B}}$ be a divisor on algebraic Hyper-K\"{a}hler
manifold N such that the linear system $|D_{\text{B}}|$ defines the
holomorphic map $\left(  \ref{Lagfibr}\right)  $ then
\begin{equation}
\mathcal{B}\text{(}P(D_{\text{B}}),P(D_{\text{B}}))=0. \label{AC2}%
\end{equation}
(See also \cite{Mat}.)
\end{lemma}

\textbf{Proof: }If I prove that
\begin{equation}
\int\limits_{\text{N}}P(D_{\text{B}})\wedge\Omega_{\text{N}}^{n-1}%
\wedge\overline{\Omega}_{\text{N}}^{n}=\int\limits_{\text{N}}P(D_{\text{B}%
})\wedge\overline{\Omega}_{\text{N}}^{n-1}\wedge\Omega_{\text{N}}^{n}=0
\label{AC20}%
\end{equation}
and
\begin{equation}
\int\limits_{\text{N}}P(D_{\text{B}})\wedge P(D_{\text{B}})\wedge
\Omega_{\text{N}}^{n-1}\wedge\overline{\Omega}_{\text{N}}^{n=1}=0 \label{AC21}%
\end{equation}
then Lemma \ref{ac-1} will follow from the definition of the
Beauville-Bogomolov form given by $\left(  \ref{Bog-B}\right)  $. Since
$P(D_{\text{B}})$ is of type $(1,1)$ and $\overline{\Omega}^{n-1}\wedge
\Omega^{n}$ is of type $(2n,2n-2),$ I can conclude $\left(  \ref{AC20}\right)
.$ Poincare duality implies that
\begin{equation}
\int\limits_{\text{N}}P(D_{\text{B}})\wedge P(D_{\text{B}})\wedge
\Omega_{\text{N}}^{n-1}\wedge\overline{\Omega}_{\text{N}}^{n=1}=\int
\limits_{\pi^{\ast}(D_{1}\cap D_{2})}\Omega_{\text{N}}^{n-1}\wedge
\overline{\Omega}_{\text{N}}^{n=1}, \label{AC21i}%
\end{equation}
where $D_{i}$ are very ample divisors on B that intersect transversely in a
subvariety of codimension two. From the fact that the fibres of the map
$\left(  \ref{Lagfibr}\right)  $ are Lagrangian, I can represent the tangent
space $T_{x,\pi^{-1}(D_{1}\cap D_{2})}$ at any non singular point $x\in
\pi^{\ast}(D_{1}\cap D_{2})$ as follows:
\begin{equation}
T_{x,\pi^{-1}(D_{1}\cap D_{2})}=V_{1}\oplus V_{2}, \label{ac14}%
\end{equation}
where $V_{1}$ is a tangent space to the fibre through the point $x$ and
$V_{2}$ is some complementary subspace in $T_{x,\pi^{-1}(D_{1}\cap D_{2})}$
that satisfies $\left(  \ref{ac14}\right)  $. Clearly
\[
\dim_{\mathbb{C}}V_{1}=\frac{1}{2}\dim_{\mathbb{C}}\text{N}=n\text{ and }%
\dim_{\mathbb{C}}V_{2}=n-1\text{.}%
\]
From the fact that $V_{1}$ is a totally isotropic $n$ dimensional subspace in
the $2n$ dimensional vector space $T_{x,\text{N}}$ with respect to the
holomorphic two form $\Omega_{\text{N}}^{2}$ I can conclude that for each
$v_{2}\in V_{2}$ we have
\[
\Omega_{\text{N}}^{2}(V_{1},v_{2})\neq0.
\]
Since the dimension of $V_{2}$ is $n-2$ then
\[
\Omega_{\text{N}}^{n-1}|_{\pi^{\ast}(D_{1}\cap D_{2})}=0.
\]
This proves formula $\left(  \ref{AC21}\right)  .$ Lemma \ref{ac-1} is proved.
$\blacksquare$

\begin{lemma}
\label{ac-3} Suppose that N can be fibered by Lagrangian manifolds. Then the
rank of $H^{2}($N,$\mathbb{Z})$ is bigger or equal to 5.
\end{lemma}

\textbf{Proof: }The definition of the Hyper-K\"{a}hler manifold implies that
\[
\dim_{\mathbb{C}}H^{2}(\text{N,}\mathbb{C})\geq3.
\]
Lemma \ref{ac-1} implies that there exists a divisor $D_{\text{B}}$ on N which
satisfies $\left(  \ref{AC2}\right)  .$ From \ the existence of $D_{\text{B}}$
that satisfies $\left(  \ref{AC2}\right)  $, the fact that the
Beauville-Bogomolov form is a non-degenerate form of signature $(3,b_{2}-3)$
and the Poincare duality I can conclude that
\[
\dim_{\mathbb{C}}H^{2}(\text{N,}\mathbb{C})\geq4.
\]
It is easy to see that since the signature of the Bogomolov-Beauville form
$\mathcal{B}$ is $(3,b_{2}-3)$ and $P(D_{\text{B}})$ is an isotropic vector$,$
then there exists $L\in H^{2}($N$,\mathbb{Z})$ such that
\begin{equation}
\mathcal{B}(L,L)>0\text{ and }\mathcal{B}(L,P(D_{\text{B}}))=0. \label{AC21a}%
\end{equation}
Then according to \cite{Hyu}, $L$ can be realized as an ample divisor on some
Hyper-K\"{a}hler manifold N.

\begin{proposition}
\label{ac-31}Suppose $L\in H^{2}($N$,\mathbb{Z})$ and it satisfies the
inequality in $\left(  \ref{AC21a}\right)  $. Suppose that N is a algebraic
Hyper-K\"{a}hler manifold with a polarization class $L$ such that
\begin{equation}
\dim_{\mathbb{C}}H^{2}(\text{N},\mathbb{Q})=4. \label{DIM}%
\end{equation}
Suppose that there exists an element $P(D_{\text{B}})\in H^{2}($%
N$,\mathbb{Z})$ such that it satisfies the equality $\left(  \ref{AC21a}%
\right)  $. Then the dimension of the coarse moduli space $\mathfrak{M}_{L}$
of polarized algebraic Hyper-K\"{a}hler manifolds is one and the Baily-Borel
compactification of the period domain
\[
\Gamma_{L}\backslash SO_{0}(2,1)/SO(2)
\]
contains a zero dimensional cusp $k_{\infty}.$
\end{proposition}

\textbf{Proof: }Let $\mathfrak{M}_{L}$ be the coarse moduli space of the
polarized algebraic Hyper-K\"{a}hler manifolds with the polarization class
$L.$ The local Torelli Theorem and the condition $\left(  \ref{DIM}\right)  $
imply that the dimension of $\mathfrak{M}_{L}$ is one. Condition $\left(
\ref{DIM}\right)  $ combined with local Torelli Theorem imply that the period
map will be a finite map
\[
p:\mathfrak{M}_{L}\rightarrow\Gamma_{L}\backslash SO_{0}(2,1)/SO(2).
\]
According to \cite{Sc}, the subgroup in $\Gamma_{L}$ that stabilizes the
totally isotropic sublattice spanned by $P(D_{\text{B}})$ that satisfies
$\left(  \ref{AC21a}\right)  $, contains a unipotent element. According to a
Theorem of Borel the existence of a unipotent element in $\Gamma_{L}$ implies
that the space $\Gamma_{L}\backslash SO_{0}(2,1)/SO(2)$ is not compact one.
This implies that the Baily-Borel compactification of
\[
\Gamma_{L}\backslash SO_{0}(2,1)/SO(2)=\mathbb{SL}_{2}(\mathbb{R})/SO(2)
\]
contains a zero dimensional cusp $k_{\infty}.$ This proves Proposition
\ref{ac-31}. $\blacksquare$

Lemma \ref{ac-3} will follow directly from the following Proposition:

\begin{proposition}
\label{ac31}Suppose that
\begin{equation}
\pi:\mathcal{N\rightarrow D} \label{ac31a}%
\end{equation}
is a family of algebraic Hyper-K\"{a}hler manifolds over the unit disk.
Suppose that the monodromy operator T$_{2}$ of the family $\left(
\ref{ac31a}\right)  $ acting on $H^{2}($N$_{t},\mathbb{Q})$ is such that
\begin{equation}
(T_{2}-id)^{2}=0 \label{Sch1}%
\end{equation}
and
\begin{equation}
T_{2}-id\neq0. \label{Sch2}%
\end{equation}
Then the Jordan normal form of the monodromy operator will contain an even
number of Jordan blocks of dimension 2.
\end{proposition}

\textbf{Proof: }The proof is based on the existence of limiting mixed Hodge
structure as proved by W. Schmidt in \cite{Sch}. In this case the weight
filtration is defined by the monodromy operator T$_{2}.$ The conditions
$\left(  \ref{Sch1}\right)  $ and $\left(  \ref{Sch2}\right)  $ imply that the
topological weight filtration is
\[
W_{0}=H^{2}(\text{N,}\mathbb{C})\supset W_{1}\supset W_{2}\supset0,
\]
where
\[
W_{1}:=\ker\ln(T)\text{ and }W_{2}:=\operatorname{Im}\ln(T).
\]
Since $W_{1}/W_{2}$ admits a Hodge structure of weight one according to
\cite{Sch}, I can conclude that
\begin{equation}
W_{1}/W_{2}=H_{W}^{1,0}\oplus\overline{H_{W}^{1,0}}. \label{Sch10}%
\end{equation}
The existence of the cusp and $\left(  \ref{Sch10}\right)  $ imply that
\begin{equation}
\dim_{\mathbb{C}}W_{2}/W_{1}=2k>0. \label{Sch11}%
\end{equation}
The relation $\left(  \ref{Sch11}\right)  $ combined with the fact that the
number of Jordan blocks of dimension 2 is equal to $\dim_{\mathbb{C}}%
W_{2}/W_{1}=2k>0$ implies Proposition \ref{ac31}. $\blacksquare$

According to Lemma \ref{ac-1}, the assumption that N can be fibered by
Lagrangian tori implies the existence of a divisor $D_{\text{B}}$ such that
its Poincare dual $P([D_{\text{B}}])\in H^{2}($N,$\mathbb{Z})$ satisfies
$\left(  \ref{AC2}\right)  $. Earlier we pointed out that this fact together
with the fact that $\Gamma_{L}$ is an arithmetic group imply the existence of
$\phi\in\Gamma_{L}$ which stabilized the isotropic vector $P([D_{\text{B}}]$
so we have%
\[
\phi(P([D_{\text{B}}])=P([D_{\text{B}}]).
\]
Let $T_{\phi}$ be the matrix of the action of the diffeomorphism $\phi$ on
$H^{2}($N,$\mathbb{Z}).$ Then $T_{\phi}$ will preserve the Bogomolov-Beauville
form. According to \cite{Sc} some power $T_{\phi}^{N_{0}}$ of $T_{\phi}$ will
be a unipotent element in the group $\Gamma_{L}\subset$ $SO_{0}(2,1).$ Since
$SO(2,1)$ is a double covering of $\mathbb{SL}_{2}(\mathbb{R})$ I can conclude
from Clemens' theory that the Jordan form of $T_{\phi}^{N_{0}}$ will contains
blocks of at most dmension one, i.e.
\[
\left(
\begin{array}
[c]{cc}%
1 & 1\\
0 & 1
\end{array}
\right)  .
\]
Let $k_{\infty}$ be the cusp in the Baily-Borel compactification
\[
\mathcal{\hat{M}}_{L}:=\overline{\Gamma_{L}\backslash SO_{0}(2,1)/SO_{0}(2)}%
\]
of $\Gamma_{L}\backslash SO_{0}(2,1)/SO(2)$ that corresponds to $T_{\phi}%
\in\Gamma_{L}.$ Repeating the arguments of Theorem \ref{j} which is proved in
the next \textbf{Section,} I can conclude that there exists a family of
polarized algebraic Hyper-K\"{a}hler manifolds%
\begin{equation}
\pi:\mathcal{N}^{\ast}\mathcal{\rightarrow D}^{\ast} \label{fhk}%
\end{equation}
over the unit disk such that the period map%
\[
p^{\ast}:\mathcal{D}^{\ast}\rightarrow\Gamma_{L}\backslash SO_{0}(2,1)/SO(2)
\]
can be prolonged to a holomorphic map
\[
p:\mathcal{D}\rightarrow\overline{\Gamma_{L}\backslash SO_{0}(2,1)/SO(2)}%
\]
of the disk $\mathcal{D}$ into the Baily-Borel compactification
\[
\mathcal{\hat{M}}_{L}:=\overline{\Gamma_{L}\backslash SO_{0}(2,1)/SO(2)}%
\]
of $\Gamma_{L}\backslash SO_{0}(2,1)/SO(2)$ where $p(0)=\kappa_{\infty},$ and
$\kappa_{\infty}$ is the cusp corresponding to $T_{\phi}\in\Gamma_{L}.$
$\left(  \ref{Sch10}\right)  $ implies that the operator $T_{\phi}$ will
contain at least two Jordan blocks of dimension 2. From here I can conclude
that that the family $\left(  \ref{fhk}\right)  $ satisfies the conditions
$\left(  \ref{Sch1}\right)  $ and $\left(  \ref{Sch2}\right)  $. From here I
can conclude that $H^{2}($N,$\mathbb{Z})$ contains two isotropic vectors with
respect to the Beauville-Bogomolov form $\mathcal{B}.$ This implies that the
rank of $H^{2}($N,$\mathbb{Z})$ is bigger or equal to 5. Proposition
\ref{ac31} is proved. $\blacksquare$

Proposition \ref{ac31} implies Lemma \ref{ac-3}. $\blacksquare$

Lemma \ref{ac-3} implies Theorem \ref{Top}. $\blacksquare$

\begin{remark}
It is not difficult to see that Hironaka's Theorem about resolution of
singularities shows that the proof of Theorem \ref{Top} works also without the
assumption that the base of the family $\left(  \ref{Lagfibr}\right)  $ is
smooth. In another paper I will prove that B is always smooth. It seems that
it is reasonable to expect that B is isomorphic to $\mathbb{CP}^{n}.$
\end{remark}

\begin{theorem}
\label{main}Suppose that N is a Hyper-K\"{a}hler manifold which admits a
Lagrangian fibration. Then the moduli space of N admits a large radius limit.
\end{theorem}

\textbf{Proof:} Theorem \ref{Top} implies that the rank of $H^{2}%
($N,$\mathbb{Z})$ is greater than $4.$ Now we can apply Theorem \ref{Cond} to
conclude Theorem \ref{main}. $\blacksquare$

\section{Baily-Borel Compactification of the Period Domain and the Geometric
Meaning of the Large Radius Limit\label{GMLR}}

\begin{theorem}
\label{j}Suppose that $\dim_{\mathbb{R}}H^{2}($N,$\mathbb{R)\geq}5,$ then
there exists a family of polarized algebraic Hyper-K\"{a}hler varieties%
\begin{equation}
\pi:\mathcal{X\rightarrow}D \label{fam}%
\end{equation}
over the unit disk $D$ with only singular fibre $\pi^{-1}(0)=$N$_{0},$ such
that the monodromy operator $T$ of the family $\left(  \ref{fam}\right)  $
acting on $H_{4n-2}($N$_{t},\mathbb{Z})$ contains a Jordan block of dimension
3, where N$_{t}=\pi^{-1}(t)$ and $t\neq0$. This means that there exists a
unique primitive cycle $\delta_{0}\in H_{4n-2}($N$_{t},\mathbb{Z})$ up to a
sign and cycles $\delta_{1}$ and $\delta_{2}$ such that
\begin{equation}
T(\delta_{0})=\delta_{0},\text{ }T(\delta_{1})=\delta_{0}+\delta_{1}\text{ and
}T(\delta_{2})=\delta_{0}+\delta_{1}+\delta_{2}. \label{BB101}%
\end{equation}

\end{theorem}

\textbf{Proof: }I will outline the proof of the Theorem first. I will prove
that there exists a polarization class $L$ such that the group $\Gamma_{L}$
contains a unipotent element $T$ that satisfies $\left(  \ref{BB101}\right)
$. Let $\overline{\mathfrak{h}_{2,b_{2}-3}/\Gamma_{L}}$ be the Baily-Borel
compactification of the locally symmetric space $\mathfrak{h}_{2,b_{2}%
-3}/\Gamma_{L}.$ The existence of the unipotent element $T$ that satisfies
$\left(  \ref{BB101}\right)  $ is equivalent to the existence of a zero
dimensional cusp $k_{\infty},$ where
\begin{equation}
k_{\infty}\in\mathfrak{D}_{\infty}:=\overline{\mathfrak{h}_{2,b_{2}-3}%
/\Gamma_{L}}\text{ }\ominus\mathfrak{h}_{2,b_{2}-3}/\Gamma_{L}. \label{Un1}%
\end{equation}
Let
\[
\mathcal{D\subset}\mathbb{C}%
\]
be a disk such that
\[
\mathcal{D}\subset\overline{\mathfrak{h}_{2,b_{2}-3}/\Gamma_{L}},\text{
}k_{\infty}\in\mathcal{D}\text{ and }\mathcal{D}^{\ast}:=\mathcal{D\ominus
}k_{\infty}\subset\mathfrak{h}_{2,b_{2}-3}/\Gamma_{L}.
\]
Then I will construct a family of polarized Hyper-K\"{a}hler manifolds over
the unit punctured disk $\mathcal{D}_{1}^{\ast}$ such that the image of base
of the family
\[
\pi^{\ast}:\mathcal{X}^{\ast}\mathcal{\rightarrow D}_{1}^{\ast}=\mathcal{D}%
_{1}\ominus0
\]
under the period map $\rho$ is exactly $\mathcal{D}^{\ast}.$

\begin{lemma}
\label{U} There exists a class of cohomology $L\in H^{2}($N,$\mathbb{Z}$) such
that $L$ can be realized as a polarization class on some algebraic
Hyper-K\"{a}hler manifold N and the group $\Gamma_{L}$ defined by Definition
\ref{MSG} contains a unipotent element $T$ that satisfies $\left(
\ref{BB101}\right)  .$
\end{lemma}

\textbf{Proof: }The proof of Lemma \ref{U} is based on Theorem \ref{Me}.
Applying Theorem \ref{Me} to the Euclidean lattice $H^{2}($N,$\mathbb{Z})/Tor$
with respect to Bogomolov-Beauville form, I can find a vector $\delta_{0}%
\in\Lambda$ such that
\begin{equation}
\mathcal{B}\left(  \delta_{0},\delta_{0}\right)  =0. \label{eq0}%
\end{equation}
According to \cite{Hyu}, if $L\in H^{2}($N,$\mathbb{Z})$ and
\begin{equation}
\mathcal{B}\left(  L,L\right)  >0 \label{eq4}%
\end{equation}
then $\mathcal{P}_{L}$ as defined in $\left(  \ref{eq3}\right)  $ will contain
an everywhere dense and open subset $\mathfrak{P}_{L}$ whose points will
correspond under the period map to marked algebraic polarized Hyper-K\"{a}hler
manifolds. Surjectivity of the period map allows us to conclude that
$\mathfrak{P}_{L}$ is not empty. See \cite{Hyu}. Since $\delta_{0}$ is an
isotropic vector, rk $H^{2}($N,$\mathbb{Z})>4$ and the signature of the
Bogomolov-Bauville form $\mathcal{B}$ is $(3,b_{2}-3)$ then one can always
find $L\in H^{2}($N,$\mathbb{Z})$ such it satisfies $\left(  \ref{eq4}\right)
$ and
\begin{equation}
\mathcal{B}\left(  \delta_{0},L\right)  =0. \label{eq5}%
\end{equation}
Suppose that it satisfies $\left(  \ref{eq4}\right)  $ and $\left(
\ref{eq5}\right)  .$ According to \cite{Sc} the subgroup in $\Gamma_{L}$ that
stabilizes the totally isotropic sublattice spanned by $\delta_{0}$, contains
an element $T_{\phi}$ such that its Jordan normal form contains a Jordan block
of dimension 3. So one can find an element $T_{\phi}\in\Gamma_{L}$ such that
there exist $\delta_{0},$ $\delta_{1}$ $and$ $\delta_{2}\in H^{2}%
($N,$\mathbb{Q})$ where $\delta_{0}$ is a primitive vector that satisfies
$\left(  \ref{eq0}\right)  $ and $T(\delta_{0})=\delta_{0},T(\delta
_{1})=\delta_{0}+\delta_{1}$ and $T(\delta_{2})=\delta_{0}+\delta_{1}%
+\delta_{2}.$ Lemma \ref{U} is proved. $\blacksquare$

\begin{lemma}
\label{U2} Suppose that
\[
b_{2}>4.
\]
Then the Baily-Borel compactification $\overline{\mathfrak{h}_{2,b_{2}%
-3}/\Gamma_{L}}$ of $\mathfrak{h}_{2,b_{2}-3}/\Gamma_{L}$ contains a zero
dimensional cusp.
\end{lemma}

\textbf{Proof: }According to Lemma \ref{U}, there exists a polarization class
$L\in H^{2}($N,$\mathbb{Z}$) such that the arithmetic group $\Gamma_{L}$
defined in Definition \ref{MSG}, contains a unipotent element. A Theorem of
Borel states that $\mathfrak{h}_{2,b_{2}-3}/\Gamma_{L}$\ is not a compact if
$\Gamma_{L}$ contains a unipotent element. See \cite{Bor}. The Baily-Borel
theory of compactifications of symmetric domains of type IV states that the
boundary $\mathcal{D}_{\infty}$ defined by $\left(  \ref{Un1}\right)  $ has
complex dimension 1. Suppose that $\mathfrak{h}_{2,b_{2}-3}/\Gamma_{L}$ is not
compact and
\[
\dim_{\mathbb{C}}\mathfrak{h}_{2,b_{2}-3}/\Gamma_{L}>1.
\]
Then $\mathcal{D}_{\infty}$ contains finite number of one dimensional and zero
dimensional cusps. See\ $\cite{BB}.$ According to \cite{BB}, the one
dimensional cusps correspond to the different orbits under the action of
$\Gamma_{L}$ of totally isotropic sublattices of rank 2 in $H^{2}%
($N,$\mathbb{Z})$. The set of zero dimensional cusps $\{k_{\infty}\}$
corresponds to the set of orbits of totally isotropic sublattices
$E_{k_{\infty}}$ of rank 1. Each such orbit that corresponds to a zero
dimensional cusp is generated \ by a nonzero primitive vector $\delta_{0}$
such that
\begin{equation}
\mathcal{B}(\delta_{0},L)=0\text{ and }\mathcal{B}(\delta_{0},\delta_{0})=0.
\label{BBII}%
\end{equation}
I already proved the existence of a vector
\[
\delta_{0}\in H_{L}^{2}(\text{N},\mathbb{Z}\text{)}:=\left\{  v\in
H^{2}(\text{N},\mathbb{Z}\text{)}|\text{ }\mathcal{B}(v,L)=0\right\}
\]
that satisfies $\left(  \ref{BBII}\right)  $. Lemma \ref{U2} is proved.
$\blacksquare$

It was already established the existence of the zero dimensional cusps if
$rkH_{2}($N,$\mathbb{Z})>4.$ Next I will show that the existence of a zero
dimensional cusp implies the existence of the family%
\begin{equation}
\pi:\mathcal{X\rightarrow}D \label{BBIII}%
\end{equation}
as stated in Theorem \ref{j}. First I will construct a family of polarized
Hyper-K\"{a}hler varieties $\left(  \ref{BBIII}\right)  $ over the punctured
unit disk $D^{\ast},$ such that the monodromy operator $T$ of the family
acting on $H_{4n-2}($N$_{t},\mathbb{Z})$ contains a Jordan block of dimension
3, where N$_{t}=\pi^{1}(t)$ and $t\neq0.$ The construction of the family
$\left(  \ref{BBIII}\right)  $ is based on a Theorem of Borel about the
existence of a subgroup of a finite index in an arithmetic group which acts
freely on a symmetric space. See \cite{Rug}.

\begin{lemma}
\label{U3}Let $\mathcal{D}$ be a disk in $\overline{\mathfrak{h}_{2,b_{2}%
-3}/\psi(\Gamma_{L})}$ such that
\[
\mathcal{D}\text{ }\ominus\text{ }0=\mathcal{D}^{\ast}=\mathcal{D}%
\cap\mathfrak{h}_{2,b_{2}-3}/\Gamma_{L}\text{ and }\mathcal{D}\cap
\mathfrak{D}_{\infty}=k_{\infty}.
\]
Then there exists a family of polarized algebraic Hyper-K\"{a}hler manifolds
$\mathcal{X}^{\ast}\rightarrow D^{\ast}$ such that $\rho(D^{\ast}%
)=\mathcal{D}^{\ast}.$ ($\rho$ is the period map.)
\end{lemma}

\textbf{Proof: }The proof of Lemma \ref{U3} is based on the following
construction of the moduli space of algebraic polarized manifolds with zero
canonical class which can be found in \cite{LTYZI}. First we proved in
\cite{LTYZI} the existence of Teichm\"{u}ller space $\mathcal{T}_{L}($N) of
polarized of algebraic polarized manifolds with zero canonical class. Then we
prove that there exists a family of marked polarized algebraic manifolds with
canonical zero
\begin{equation}
\mathcal{X}_{L}\rightarrow\mathcal{T}_{L}(\text{N}) \label{grif}%
\end{equation}
over the Teichm\"{u}ller space $\mathcal{T}_{L}($N) which is defined uniquely
up to the action of a finite group $G$ of automorphisms of N which preserve
the polarization class $L$ and acts as identity on the middle cohomology of N.

One of the main results proved in \cite{Hyu} is the global Torelli Theorem for
polarized algebraic Hyper-K\"{a}hler manifolds. The global Torelli Theorem
implies that the moduli space of polarized Hyper-K\"{a}hlerian manifold
$\mathfrak{M}_{L}($N) is biholomorphic to an open and everywhere dense subset
in $\Gamma_{L}\backslash\mathfrak{h}_{2,b_{2}-3}$ where
\[
\mathfrak{h}_{2,b_{2}-3}=SO(2,b_{2}-3)/SO(2)\times SO(b_{2}-3).
\]
It is proved in \cite{Rug} the existence of a subgroup $\Gamma_{L}%
^{^{^{\prime}}}$ of finite index in $\Gamma_{L}$ such that
\begin{equation}
\mathfrak{h}_{2,b_{2}-3}/\Gamma_{L}^{^{\prime}}=\mathcal{M}_{L}(\text{N})
\label{grifa}%
\end{equation}
is a non singular variety. We proved in \cite{LTYZI} that the moduli space of
polarized algebraic Hyper-K\"{a}hler manifolds
\begin{equation}
\mathfrak{M}_{L}(\text{N})=\Gamma_{L}^{^{\prime}}\backslash\mathcal{T}%
_{L}(\text{N)} \label{grif0}%
\end{equation}
is the quotient of the Teichm\"{u}ller space by the group $\Gamma_{L}$. Since
the subgroup $\Gamma_{L}^{^{"}}$ acts without fixed points on the
Teichm\"{u}ller space $\mathcal{T}_{L}$(N), then the family $\left(
\ref{grif}\right)  $ defines a family of polarized algebraic Hyper-K\"{a}hler
manifolds%
\begin{equation}
\mathcal{Y}_{L}\rightarrow\mathcal{M}_{L}(\text{N)}=\Gamma_{L}^{^{\prime}%
}\backslash\mathcal{T}_{L}(\text{N).} \label{grif1}%
\end{equation}
Let me recall that $\mathcal{M}_{L}($N) is a finite cover of the course moduli
space of polarized Hyper-K\"{a}hler manifolds $\mathfrak{M}_{L}($N). According
to the result of Viehweg, $\mathfrak{M}_{L}($N) is a quasi-projective
manifold. So $\mathcal{M}_{L}($N) is a quasi-projective manifold too. Let
$\overline{\mathcal{M}_{L}}$ be some projective smooth variety such that
\[
\overline{\mathcal{M}_{L}(\text{N)}}\text{ }\ominus\text{ }\mathcal{M}%
_{L}(\text{N)}=\mathfrak{D}_{\infty}%
\]
is a divisor with normal crossings. $\overline{\mathcal{M}_{L}(\text{N)}}$
exists by Hironaka's Theorem of resolution of singularities. According to a
Theorem of Borel in \cite{Bo}, there exists a projective map
\[
\overline{\rho}:\overline{\mathcal{M}_{L}(\text{N)}}\rightarrow\overline
{\mathfrak{h}_{2,b_{2}-3}/\Gamma_{L}}%
\]
such that the map
\[
\rho:\mathcal{M}_{L}(\text{N)}\rightarrow\mathfrak{h}_{2,b_{2}-3}/\Gamma_{L}%
\]
is a finite map. The local Torelli Theorem implies that $\overline{\rho}$ is a
surjective map. Then the restriction of the family $\left(  \ref{grif1}%
\right)  $ over the lift $D^{\ast}$
\[
D^{\ast}\subset\mathcal{T}_{L}(\text{N})/\Gamma_{L}^{^{\prime}}=\mathcal{M}%
_{L}\text{(N)}%
\]
\ of $\ \mathcal{D}^{\ast}\subset\mathfrak{h}_{2,b_{2}-3}/\Gamma_{L}$ defines
a family
\begin{equation}
\pi^{\ast}:\mathcal{X}^{\ast}\rightarrow D^{\ast} \label{grif2}%
\end{equation}
of polarized algebraic Hyper-K\"{a}hler varieties. From Griffiths' results
proved in \cite{Grif}, it follows that the restriction of the period map
$\rho^{^{\prime}}$ to $D^{\ast}$ can be prolonged to a map%
\begin{equation}
\overline{\rho^{^{\prime}}}:D\rightarrow\overline{\mathfrak{h}_{2,b_{2}%
-3}/\psi(\Gamma_{L})}, \label{grif2a}%
\end{equation}
where $\overline{\mathfrak{h}_{2,b_{2}-3}/\psi(\Gamma_{L})}$ is the
Baily-Borel compactification of $\mathfrak{h}_{2,b_{2}-3}/\psi(\Gamma_{L}%
).$See \cite{LTY}. Lemma \ref{U3} is proved. $\blacksquare$

\begin{remark}
It is easy to see that $\left(  \ref{grif2a}\right)  $ implies that the family
$\left(  \ref{grif2}\right)  $ can be compactified to a family $\pi
:\mathcal{X\rightarrow}D$ of algebraic manifolds over the unit disk$.$
\end{remark}

It is a well known fact that the monodromy operator $T_{\phi}$ acting on
$H_{4n-2}($N,$\mathbb{Z})$ or $H^{2}($N$,\mathbb{Z})$ of the family
\[
\pi^{\ast}:\mathcal{X}^{\ast}\mathcal{\rightarrow}D^{\ast}%
\]
constructed in Lemma \ref{U3} will be the operator $T_{\phi}$ induced by the
element $\phi\in\Gamma$(N) is induced by the geometric monodromy $\phi$ which
is an element of the mapping class group $\Gamma$(N). $T_{\phi}$ will contain
a Jordan block of dimension 3. Clemens theory of monodromy implies that there
exist $4n-2$ cycles $\delta_{0},$ $\delta_{1}$ and $\delta_{2}$ such that it
satisfies $\left(  \ref{BB101}\right)  $. (For details about Clemens' theory
see \cite{Cl}.) This follows from the way I defined $\pi^{\ast}:\mathcal{X}%
^{\ast}\mathcal{\rightarrow}D^{\ast}.$ Theorem \ref{j} is proved.
$\blacksquare$

\begin{corollary}
Definition \ref{LRL} is equivalent to Definition \ref{LRLI}.
\end{corollary}

\begin{definition}
\label{priminv}Let $\delta_{0}$ be a primitive $4n-2$ cycle that satisfies
$\left(  \ref{BB101}\right)  $. I will call $\delta_{0}$ a primitive invariant
cycle associated with the monodromy operator $T.$ ($\delta_{0}$ is defined up
to a sign.)
\end{definition}

\begin{corollary}
\label{BB1}Suppose that N is a Hyper-K\"{a}hler manifold which satisfies the
conditions of Theorem \ref{j}. Let $\delta$ be a 4n-2 dimensional non zero
cycle on N. Suppose that the Poincare dual $P(\delta)$ of $\delta$ satisfies
\begin{equation}
\mathcal{B}(P(\delta),P(\delta)\mathcal{)}=0. \label{AC33}%
\end{equation}
Suppose that $\delta=m\delta_{0},$ where $\delta_{0}$ is a primitive element
in $H_{4n-2}($N,$\mathbb{Z}).$ Then the condition $\left(  \ref{AC33}\right)
$ is necessary and sufficient for the existence of a zero dimensional cusp
$k_{\infty}$ of the Baily-Borel compactification of the period domain and
$k_{\infty}$ corresponds to the vanishing invariant cycle $\delta_{0}.$
\end{corollary}

Clemens theory of Picard-Lefschetz transformations implies:

\begin{corollary}
\label{BB12}Suppose that $\delta$ is a 4n-2 dimensional non zero cycle on N.
Suppose that it satisfies $\left(  \ref{AC33}\right)  $. Then%
\[
\underset{k}{\underbrace{\left\langle \delta,...,\delta\right\rangle }}=0
\]
for $k>n$ and
\[
\underset{n}{\underbrace{\left\langle \delta,...,\delta\right\rangle }}\neq0
\]
See \cite{Cl} or \cite{LTY}.
\end{corollary}

\section{The Structure of the Lagrangian Fibrations of Hyper-K\"{a}hler
Manifolds}

Matsushita proved in \cite{Mat} that if N is an algebraic Hyper-K\"{a}hler
manifold that admits a fibration then it is a Lagrangian holomorphic fibration
and the base has the same Hodge numbers as the projective space $\mathbb{CP}%
^{n}$.

\begin{theorem}
\label{ac}Let N be a Hyper-K\"{a}hler manifold such that it admits a
Lagrangian fibration
\begin{equation}
\pi:\text{N}\rightarrow\text{B.} \label{AC0}%
\end{equation}
Then the generic fibre of $\left(  \ref{AC0}\right)  $ is an n dimensional
complex tori. According to Theorem \ref{Top} B is a projective algebraic
variety. Suppose that $D_{\text{B}}$ is an ample divisor on B then $\pi^{\ast
}(D_{\text{B}})$ as a $4n-2$ cycle is homological to a vanishing invariant
$4n-2$ cycle of the action of some element $\phi\in\Gamma$(N) induced on
$H_{4n-2}($N,$\mathbb{Z})$. Moreover the homology class of the fibre is
homological to the vanishing invariant cycle of the unique monodromy operator
T of maximal unipotent rank which acts on $H_{2n}($N,$\mathbb{Z}).$
\end{theorem}

\textbf{Proof: }Theorem \ref{Top} implies that if N admits a Lagrangian
fibration then the second homology group will be of rank greater or equal to
five and the base B will be a projective algebraic variety. So I can apply
Theorem \ref{j} to conclude the existence of a maximal unipotent element
$\phi$ in the mapping class group which has a Jordan block of rank 3 such that
the induced action $T_{\phi}$ of $\phi$ on $H^{2}($N,$\mathbb{Z})$ has the
following properties$:$%
\[
T_{\phi}(\mathcal{P}(\pi^{\ast}(D_{\text{B}}))=\mathcal{P}(\pi^{\ast
}(D_{\text{B}})
\]
and there exist elements $\delta_{1}$ and $\delta_{2}$ of $H^{2}%
($N,$\mathbb{Z})$ such that
\[
T_{\phi}(\delta_{1})=\delta_{1}+\pi^{\ast}(D_{\text{B}})\text{ and }T_{\phi
}(\delta_{2})=\delta_{2}+\delta_{1}.
\]

In order to prove Theorem \ref{ac} I need to check that the $4n-2$ cycle
$\pi^{\ast}(D_{\text{B}})$ as defined in Theorem \ref{ac}, satisfies the
condition $\left(  \ref{AC33}\right)  $ of Cor. \ref{BB1}. The condition
$\left(  \ref{AC33}\right)  $\ was already proved in Lemma \ref{ac-1}. This
implies that $\pi^{\ast}(D_{\text{B}})$ as a $4n-2$ cycle is homological to a
vanishing invariant $4n-2$ cycle of some unipotent element of the mapping
class group of N acting on $H_{4n-2}($N,$\mathbb{Z})$ in the family $\left(
\ref{fam}\right)  $ constructed in Theorem \ref{j}.

From the fact that $D_{\text{B}}$ is a very ample divisor on B, I can conclude
that I can find $n$ divisors $D_{1},...,D_{n}$ in the linear system
$|D_{\text{B}}|$ such that they intersect each other transversely in $k$
different points and these $k$ points are not in the discriminant locus of the
fibration $\left(  \ref{AC0}\right)  .$ Let
\[
\left\langle \lbrack\pi^{\ast}(D_{1})],...,[\pi^{\ast}(D_{n})]\right\rangle
=[C_{0}],
\]
be a $n$ cycle which is homological to the n times self intersection of the
vanishing invariant $4n-2$ cycle $[\pi^{\ast}(D_{1})]$. From Lemma\ \ref{031}
and Theorem \ref{Ver} it follows that that $[C_{0}]$ is homological to the
vanishing invariant cycle of the monodromy operator $T=S^{n}(T_{\phi})$ which
is of maximal unipotent rank and acts on $H_{2n}($N,$\mathbb{Z}).$ This fact
implies that the homology class of the of the fibration $\left(
\ref{AC0}\right)  $ is homological to the vanishing invariant cycle of the
monodromy operator $T_{\phi}$. In order to finish the proof of Theorem
\ref{ac} I need to prove that the generic fibre of the fibration $\left(
\ref{AC0}\right)  $ is a complex tori.

\begin{lemma}
\label{ac0}The generic fibre
\[
\pi^{-1}(t)=F_{t}%
\]
is a tori. (See also \cite{Mat1}.)
\end{lemma}

\textbf{Proof: }Since $F_{t}$ is a generic fibre, it means that $F_{t}$ is a
non-singular complex submanifold. Notice that the normal bundle $\mathcal{N}%
_{F_{t}/\text{N}}$ is the trivial bundle. From the fact that $F_{t}$ is a
Lagrangian submanifold in N, i.e.
\[
\Omega_{\text{N}}|_{F_{t}}\equiv0
\]
implies that the cotangent bundle $\Omega_{F_{t}}^{1}$ to $F_{t}$ is
isomorphic to the normal bundle $\mathcal{N}_{F_{t}/\text{N}}.$ So this
implies that $F_{t}$ is a tori. $\blacksquare$

Theorem \ref{ac} is proved. $\blacksquare$

\section{Riemann-Roch-Hirzebruch Theorem}

Theorem \ref{j42} can be derived from the results in \cite{Mat1} and
\cite{Hyu}. A different approach in the proof of $\left(  \ref{Hirz}\right)  $
is used in this article.

\begin{theorem}
\label{j42}Let N be a Hyper-K\"{a}hler manifold such that $D$ is an
irreducible divisor that is homological to a vanishing invariant $4n-2$ cycle
$\delta_{0}$ described in Theorem \ref{j}. Then the following formula is true
\begin{equation}
\chi(\mathcal{O}_{\text{N}}(D))=\sum_{i=0}^{4}(-1)^{i}\dim_{\mathbb{C}}%
H^{i}(\text{N},\mathcal{O}_{\text{N}}(D)=n+1. \label{Hirz}%
\end{equation}

\end{theorem}

\bigskip\textbf{Proof: }The proof of Theorem \ref{j42} is based on the
Hirzebruch-Riemann-Roch Theorem:
\[
\chi(\mathcal{O}_{\text{N}}(D))=\int\limits_{\text{N}}Td(T_{\text{N}%
})ch(\mathcal{O}_{\text{N}}(D)).
\]
(See \cite{Hir}.) $\chi(\mathcal{O}_{\text{N}}(D))$ means the Euler
characteristics of the line bundle $\mathcal{O}_{\text{N}}(D),$
$Td(T_{\text{N}}$) means the Todd class of the tangent bundle of the
Hyper-K\"{a}hler manifold N and $ch(\mathcal{O}_{\text{N}}(D))$ means the
Chern class of the line bundle $\mathcal{O}_{\text{N}}(D).$ Theorem \ref{j42}
will follow directly from the following Lemma:

\begin{lemma}
\label{j421}I have the following formula:
\[
\int\limits_{\text{N}}Td(T_{\text{N}})Ch(\mathcal{O}_{\text{N}}(D))=\int
\limits_{\text{N}}Td(T_{\text{N}})_{(2n,2n)}=n+1.
\]

\end{lemma}

\textbf{Proof: }Cor. \ref{BB12} implies that
\begin{equation}
c_{1}(\mathcal{O}_{\text{N}}(D))^{n+1}=...=c_{1}(\mathcal{O}_{\text{N}%
}(D))^{2n}=0. \label{Ri}%
\end{equation}
The formula $\left(  \ref{Ri}\right)  $ implies
\begin{equation}
\left(  Td(\text{N})Ch(\mathcal{O}_{\text{N}}(D))\right)  _{(2n,2n)}%
=\sum_{k\leq n}\frac{1}{k!}Td(N)_{(2n-k,2n-k)}c_{1}(\mathcal{O}_{\text{N}%
}(D))^{k}. \label{Ria}%
\end{equation}
I will show that formula $\left(  \ref{Ria}\right)  $ implies:
\begin{equation}
\int\limits_{\text{N}}\left(  Td(\text{N})Ch(\mathcal{O}_{\text{N}%
}(D))\right)  _{(2n,2n)}=\int\limits_{\text{N}}Td(N)_{(2n,2n)}. \label{Rib}%
\end{equation}
Formula $\left(  \ref{Ria}\right)  $ implies formula $\left(  \ref{Rib}%
\right)  ,$ if the following relation
\begin{equation}
\int\limits_{\text{N}}Td(N)_{(2n-k,2n-k)}\left(  c_{1}(\mathcal{O}_{\text{N}%
}(D))\right)  ^{k}=0 \label{Ric}%
\end{equation}
is true for $0<k\leq n.$ Let
\[
\delta_{k}(t):=\underset{k}{\underbrace{\left\langle \delta_{0}(t),...,\delta
_{0}(t)\right\rangle }}%
\]
for $k\leq n.$ The fact that the Todd class is a topological invariant implies
that $Td($N)$_{(2n-k,2n-k)}$ will be always of type $(2n-k,2n-k).$ Applying
Poincare duality and the fact that N is obtained from a family $\{$N$_{t}%
$\}$\rightarrow\mathcal{D}$ over the unit disk $\mathcal{D}$ by isometric
deformation, I deduce from
\begin{equation}
P(\delta_{k}(t))=c_{1}(\mathcal{O}_{\text{N}}(D))^{k}. \label{Ricci}%
\end{equation}
Combining $\left(  \ref{Ricci}\right)  $ with the fact that $Td($%
N)$_{(2n-k,2n-k)}$ will be always of type $(2n-k,2n-k),$ I obtain for $0<k\leq
n$ the following formula:
\begin{equation}
\int\limits_{\text{N}}Td(\text{N)}_{(2n-k,2n-k)}c_{1}(\mathcal{O}_{\text{N}%
}(D))^{k}=\int\limits_{\delta_{k}(t)}Td(\text{N)}_{(2n-k,2n-k)}. \label{6}%
\end{equation}
Then $\left(  \ref{6}\right)  $ implies that the\ integral in the right-hand
side of $\left(  \ref{6}\right)  $ is a topological invariant. So
\[
\int\limits_{\delta_{k}(t)}Td(\text{N)}_{(2n-k,2n-k)}=const
\]
for all $t\in D.$ In order to compute this constant one needs to compute
\[
\underset{t\rightarrow0}{\lim}\int\limits_{\delta_{k}(t)}Td(\text{N)}%
_{(2n-k,2n-k)}.
\]
From Clemens theory of vanishing cycles it follows that if $k=n$ then
\[
\underset{t\rightarrow0}{\lim}\text{ }\delta_{k}(t)=point.
\]
For $k\leq n+1$ then $\dim\delta_{k}(t)=4n-2k$
\[
\underset{t\rightarrow0}{\lim}\text{ }\delta_{k}(t)=\gamma_{0}%
\]
where $\gamma_{0}$ is cycle on the central fibre N$_{0}$ of the of the family
described above and it has dimension less then $4n-2k.$ This implies that
\begin{equation}
\int\limits_{\text{N}}Td(\text{N)}_{(2n-k,2n-k)}c_{1}(\mathcal{O}_{\text{N}%
}(D))^{k}=\underset{t\rightarrow0}{\lim}\text{ }\int\limits_{\delta_{k}%
(t)}Td(\text{N)}_{(2n-k,2n-k)}=0. \label{Rid}%
\end{equation}
Formula $\left(  \ref{Rid}\right)  $ I obtain that
\begin{equation}
\int\limits_{\text{N}}Td(N)Ch(\mathcal{O}_{\text{N}}(D))=\int\limits_{\text{N}%
}Td(T_{\text{N}})_{(2n,2n)}=\chi(\mathcal{O}_{\text{N}}). \label{Rif}%
\end{equation}
It is known that since N is a Hyper-K\"{a}hler manifold then $\chi
(\mathcal{O}_{\text{N}})=n+1.$ This fact together with $\left(  \ref{Rif}%
\right)  $ imply that
\begin{equation}
\int_{\text{N}}Td(T_{\text{N}})_{(4,4)}=\chi(\mathcal{O}_{\text{N}})=n+1.
\label{5}%
\end{equation}
Lemma \ref{j421} is proved. $\blacksquare$\textbf{\ }

As it was pointed out already Lemma \ref{j421} implies Theorem \ref{j42}.
$\blacksquare$

\section{Final Remarks and Some Conjectures}

It is easy to prove the set of all Hyper-K\"{a}hler manifolds whose second
homology group has a rank greater than 4 and whose points correspond to
Hyper-K\"{a}hler manifolds whose Neron Severi group is generated by an
isotropic vector with respect to Beauville-Bogomolov, form an everywhere dense
subset in the moduli space of marked Hyper-K\"{a}hler manifolds. I know that
for those Hyper-K\"{a}hler manifolds I can construct a holomorphic line bundle
whose first Chern class is the isotropic element of $H^{2}($N,$\mathbb{Z}).$
The main problem of the SYZ conjecture for Hyper-K\"{a}hler manifolds is to
construct holomorphic sections of this line bundle. It seems to me that
Clemens' theory suggests that on this line bundle one can find a metric whose
Chern form has the following property; on a "large" open part of the
Hyper-K\"{a}hler manifold the Hermitian form associated with the Chern form
has exactly $n=\frac{1}{2}\dim_{\mathbb{C}}$N positive eigen values and the
rest are zero. It is natural to study if the Witten Holomorphic Morse
technique developed by J.-P. Demailly, could be used to prove the following conjecture:

\begin{conjecture}
Suppose that N is a Hyper-K\"{a}hler manifold whose Neron-Severi group has
rank one and it is generated by an isotropic element of $H^{2}($%
N,$\mathbb{Z})$ with respect to the field of meromorphic functions on N has a
transcendence degree $n=\frac{1}{2}\dim_{\mathbb{C}}$N.
\end{conjecture}

It seems that it is reasonable to suspect that there do not exist
Hyper-K\"{a}hler manifolds whose field of meromorphic functions has positive
transcendental degree less than $n.$ If I know that the transcendental degree
of the field of meromorphic function on N is equal to $n=\frac{1}{2}%
\dim_{\mathbb{C}}$N and the Neron-Severi group of N is generated by an
isotropic vector then it is not difficult to show that N can be fibered on
Lagrangian tori over $\mathbb{CP}^{n}.$

\begin{conjecture}
Suppose that N is a Hyper-K\"{a}hler manifold such that the Neron-Severi group
$H^{2}($N,$\mathbb{Z})\cap H^{1,1}($N,$\mathbb{R})$ is generated by an
isotropic vector $l.$ Then N can be fibered by Lagrangian tori over
$\mathbb{CP}^{n}$, where $\dim_{\mathbb{C}}$N=$2n.$
\end{conjecture}

It seems that the existence of the SYZ fibration for Hyper-K\"{a}hler
manifolds implies the following Conjecture:

\begin{conjecture}
There are only finite number of topological types of compact Hyper-K\"{a}hler
manifolds in each dimension.
\end{conjecture}

Another approach to the existence of SYZ fibrations on Hyper-K\"{a}hler
manifold is to try to deform the invariant cycle corresponding to the maximal
unipotent element in the mapping class group to a complex analytic manifold.
Then the deformation theory will establish the existence of SYZ fibration.

\section{Appendix}

Let $\Lambda$ be a lattice with a scalar product $\left\langle
u,v\right\rangle $ which has signature $(3,n).$ I will define a Hodge
structure of weight two on $\Lambda\mathcal{\otimes}\mathbb{C}$ as an
orthogonal decomposition of
\[
\Lambda\otimes\mathbb{C=}H^{2,0}\oplus H^{1,1}\oplus H^{0,2}%
\]
with the following properties: \textbf{1. }$\dim H^{2,0}=1.$ \textbf{2.
}$\overline{H^{2,0}}=H^{0,2}.$ \textbf{3. } $H^{1,1}=(H^{2,0}\oplus
H^{0,2})^{\perp}.$ \textbf{4. }Let $\omega\in H^{2,0}$ and $\omega\neq0.$ I
will assume that
\[
\left\langle \omega,\omega\right\rangle =0\text{ }and\text{ }\left\langle
\omega,\overline{\omega}\right\rangle >0.
\]
It is easy to see that $H^{2,0}$ defines in a unique way the Hodge
decomposition
\[
\Lambda\otimes\mathbb{C=}H^{2,0}\oplus H^{1,1}\oplus H^{0,2}.
\]
Indeed once $H^{2,0}$ is defined then $H^{0,2}$ and $H^{1,1}$ are defined
uniquely since
\[
\overline{H^{2,0}}=H^{0,2}\text{ }%
\]
and
\[
H^{1,1}=(H^{2,0}\oplus H^{0,2})^{\perp}.
\]
It is a standard fact that the group SO$_{0}$(3,n) acts transitively on the
set of all variations Hodge structures of weight two with $h^{2,0}=1$ and the
stabilizer of a fixed Hodge structure on $\Lambda\otimes\mathbb{C}$
\[
\Lambda\otimes\mathbb{C=}H^{2,0}\oplus H^{1,1}\oplus H^{0,2}%
\]
is SO($2$)$\times$SO$_{0}$($1,n$). From here I obtain that the space
$\mathcal{G}$ that parametrizes the set of all variations of Hodge structures
of weight two with $h^{2,0}=1$ is isomorphic to\
\begin{equation}
\mathcal{P}:=SO_{0}(3,n)/SO_{0}(2)\times SO(1,n). \label{Perd}%
\end{equation}

The following Theorem is standard one:

\begin{theorem}
\label{h2} The moduli space $\mathcal{P}$ of variations of Hodge structures of
weight two with $h^{2,0}=1$ is isomorphic to the open set of the quadric in
$\mathbb{P}(\Lambda\otimes\mathbb{C})$ described by
\[
\left\langle \tau,\tau\right\rangle =0\text{ and }\left\langle \tau
,\overline{\tau}\right\rangle >0.
\]

\end{theorem}

\textbf{Proof: }For the proof of Theorem \ref{h2} see \cite{To80}.
$\blacksquare$

Bogomolov proved the there are no obstructions to the deformations of complex
structures of Hyper-K\"{a}hler manifolds in \cite{B}. In \cite{B} the
Kuranishi family $\mathcal{X\rightarrow K}$(N) of Hyper-K\"{a}hler was
constructed. Since
\[
\mathcal{X}\approxeq\text{N}\times\mathcal{K}\text{(N)}%
\]
in the C$^{\infty}$ category when I marked one fibre, then all the fibres of
the family $\mathcal{X}\rightarrow\mathcal{K}$(N)(N) will be marked. Let us
fix a basis $(\gamma_{1},..,\gamma_{b_{2}})$ of $H_{2}($N$,\mathbb{Z}$)/$Tor$,
then $(\gamma_{1},..,\gamma_{b_{2}})$ will be a basis in $H_{2}($N$_{\tau
},\mathbb{Z}$)/$Tor$ for each $\tau\in\mathcal{K}$(N)(N)$.$

The period map
\[
\rho:\mathcal{K}\text{(N)}\mathcal{\rightarrow}\mathbb{P(}H^{2}(\text{N}%
_{\tau},\mathbb{Z})\otimes\mathbb{C})
\]
is defined as follows:
\begin{equation}
\rho(\tau):=\left(  ...,\int_{\gamma_{i}}\Omega_{\tau},...\right)  ,
\label{Tor}%
\end{equation}
where\textit{\ }$\Omega_{\tau}$\textit{\ }is the only holomorphic two form
defined up to a constant on N$_{\tau}:=\pi^{-1}(\tau).$

An important role will be played by the well-known local Torelli Theorem for
manifolds with canonical class zero, proved by Ph. Griffiths:

\begin{theorem}
\label{Tor1} The period map $\rho$ defined by $\left(  \ref{Tor}\right)  $ is
a local isomorphism.
\end{theorem}

\begin{criterion}
\label{r4}A cohomology class $\alpha\in H^{2}($N$,\mathbb{Z})$ is of type
$(1,1)$ on a Hyper-K\"{a}hler manifold if and only if
\begin{equation}
\int_{\text{N}}\alpha\wedge\Omega_{\text{N}}\wedge\left(  \Omega_{\text{N}%
}\wedge\overline{\Omega_{\text{N}}}\right)  ^{n-1}=0. \label{CR1}%
\end{equation}

\end{criterion}

\textbf{Proof: }The proof of Criterion \ref{r4} can be found in \cite{Bau}.
$\blacksquare$

In \cite{Bau} it is proved that
\[
\rho(\mathcal{K}\text{(N)}\mathcal{)}\subset\mathcal{P}_{1}\subset
\mathbb{P(}H^{2}(\text{N},\mathbb{Z})\otimes\mathbb{C}),
\]
where $\mathcal{P}_{1}$ is defined as one of the connected components of
$\mathcal{P}$ defined by $\left(  \ref{Perd}\right)  $, i.e.%
\[
\mathcal{P}:=SO(3,b_{2}-3)/SO(2)\times SO(1,b_{2}-3)=
\]%
\begin{equation}
\left\{  \tau\in\mathbb{P(}H^{2}(\text{N},\mathbb{Z})\otimes\mathbb{C}%
)|\mathcal{B}(\tau,\tau)\text{ \& }\mathcal{B}(\tau,\overline{\tau})\right\}
>0. \label{CR5}%
\end{equation}
Let us define%

\begin{equation}
\mathcal{P}_{L}:=\{\tau\in\mathcal{P}_{1}|\mathcal{B}(\tau,L)=0\}\text{
}and\text{ }\mathcal{K}(\text{N})_{L}:=\rho^{-1}(\mathcal{P}_{L}\cap
\rho(\mathcal{K}(\text{N}))). \label{eq3}%
\end{equation}

\begin{remark}
\label{r3}It is an elementary exercise to show that if $\mathcal{B}(L,L)>0,$
then
\[
\mathcal{P}_{L}\approxeq\mathfrak{h}_{2,b_{2}-3}=SO_{0}(2,b_{2}-3)/SO(2)\times
SO(b_{2}-3).
\]
See \cite{To80}.
\end{remark}

\begin{remark}
\label{r2} Kodaira proved the following characterization of projective
algebraic varieties; A complex compact analytic manifold N is projective iff
on N there exists a line bundle $\mathcal{L}$ and a metric $||$
$\ ||_{\mathcal{L}}^{2}$ on $\mathcal{L}$ such that $c_{1}(||\ ||_{\mathcal{L}%
}^{2})$ is point wise positive on N. (See \cite{KM}). Using this result and
Criterion \ref{r4} it is easy to see that the restriction of the Kuranishi
family
\[
\mathcal{X\rightarrow K}\text{(N)}%
\]
to
\begin{equation}
\mathcal{K}_{L}(\text{N)}=\rho^{-1}(\mathcal{P}_{L}\cap\rho(\mathcal{K(}%
\text{N)})) \label{Kur}%
\end{equation}
defines a family of algebraic and polarized Hyper-K\"{a}hler manifolds
\begin{equation}
\mathcal{X}_{L}\rightarrow\mathcal{K}_{L}(\text{N)} \label{Kur1}%
\end{equation}
for $\mathcal{K}$(N) \textquotedblright small\textquotedblright\ enouph$.$
\end{remark}

\end{document}